\PassOptionsToPackage{unicode}{hyperref}
\PassOptionsToPackage{hyphens}{url}
\documentclass[
]{article}
\usepackage{amsmath,amssymb}
\usepackage{lmodern}
\usepackage{iftex}
\usepackage{amssymb}
\usepackage{amsmath}
\usepackage{fancyhdr}
\usepackage{enumerate}
\usepackage{graphicx}
\usepackage{biblatex}
\usepackage{tikz}
\usepackage{tikz-qtree}
\usepackage{tikz-qtree,tikz-qtree-compat}
\graphicspath{{images/}}
\usepackage{mathtools}

\newcommand{\RNum}[1]{\uppercase\expandafter{\romannumeral #1\relax}}

\newcommand{\qed}{\begin{flushright}
{$\blacksquare$}
\end{flushright}
}
\ifPDFTeX
  \usepackage[T1]{fontenc}
  \usepackage[utf8]{inputenc}
  \usepackage{textcomp} 
\else 
  \usepackage{unicode-math}
  \defaultfontfeatures{Scale=MatchLowercase}
  \defaultfontfeatures[\rmfamily]{Ligatures=TeX,Scale=1}
\fi
\IfFileExists{upquote.sty}{\usepackage{upquote}}{}
\IfFileExists{microtype.sty}{
  \usepackage[]{microtype}
  \UseMicrotypeSet[protrusion]{basicmath} 
}{}
\makeatletter
\@ifundefined{KOMAClassName}{
  \IfFileExists{parskip.sty}{%
    \usepackage{parskip}
  }{
    \setlength{\parindent}{0pt}
    \setlength{\parskip}{6pt plus 2pt minus 1pt}}
}{
  \KOMAoptions{parskip=half}}
\makeatother
\usepackage{xcolor}
\IfFileExists{xurl.sty}{\usepackage{xurl}}{} 
\IfFileExists{bookmark.sty}{\usepackage{bookmark}}{\usepackage{hyperref}}
\hypersetup{
  pdftitle={Vertices of rank 1 in increasing trees and derangements},
  pdfauthor={Mario Midence Ordonez},
  hidelinks,
  pdfcreator={LaTeX via pandoc}}
\ifLuaTeX
  \usepackage{selnolig}  
\fi

\title{A bijection between vertices of rank 1 in increasing trees and derangements}
\author{Mario Midence Ordonez}
\date{}
\begin{document}
\maketitle
\begin{abstract}
We give a combinatorial proof of a recent result of B{\'o}na \cite{bona2022} by constructing a bijection from the set of all neighbors of leaves of increasing trees of size $n$ to the set of derangements of length $n$. 
\end{abstract}

{permutation, derangement, increasing tree, neighbors of leaves}

\section{1. Definitions and
Background}\label{introduction}

An \textit{increasing tree} of size $n$ is a nonplanar rooted tree with $n$ vertices in which the vertices have labels 0 through $n-1$, with the additional property that the unique downward path from the root (always labeled 0) to any vertex forms an increasing sequence. Properties and statistics of both plane and nonplanar increasing trees have been studied. It is well known that plane increasing trees are in bijection with Stirling permutations. In \cite{Janson2008a}, Janson used this correspondence to relate statistics in plane increasing trees to statistics in permutations. In \cite{Holmgren2014}, Janson and Holmgren study the probability distributions of properties of subtrees of plane increasing trees.

From here on, we will discuss the nonplanar case, so an increasing tree will be implicitly nonplanar. We first note that given an increasing tree of size $n-1$ we can construct an increasing tree of size $n+1$ by attaching the vertex $n$ as the child of any of the $n-1$ vertices. Then, if we let $T_n$ be the the number of increasing trees of size $n$, we see that $T_{n+1}=nT_{n}$. This shows that increasing trees of size $n$ are in bijection with permutations of length $n-1$. 

Before we show an explicit bijection, we need to define what a depth search walk is. A \textit{depth search walk} on an increasing tree starting at a vertex $j$ goes from $j$ to its greatest child $i$. Then we explore the subtree rooted at $i$ (recursively, using the same rules), recording each vertex. Then we return to $j$, and continue with the next greatest child of $j$, until there are no more children left. In other words, when presented with a choice, we always go to the greatest vertex that has not yet been visited.

Then the bijection between increasing trees and permutations may be defined by simply taking the depth search walk of the increasing trees, recording each vertex visited, and omitting the root 0. From this bijection we also get that leaves correspond to descents in permutations, as shown by Stanley in \cite{stanley2011enumerative}. Therefore exactly one half of the total vertices in increasing trees of size $n$ are leaves. See below the increasing tree corresponding to the permutation $4752613$.

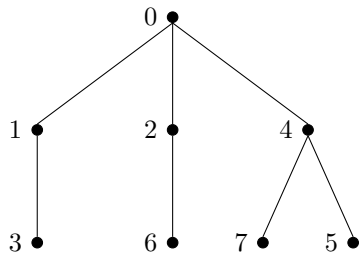
\begin{tikzpicture}
  [
    black/.style={circle, draw=black, fill=black, minimum size=4pt, inner sep=0pt},
    ublack/.style={circle, draw=black, minimum size=4pt, inner sep=0pt},
    red/.style={circle, draw=red, fill=red, minimum size=4pt, inner sep=0pt},
    ured/.style={circle, draw=red, minimum size=4pt, inner sep=0pt},
    blue/.style={circle, draw=blue, fill=blue, minimum size=4pt, inner sep=0pt},
    yellow/.style={circle, draw=yellow, fill=yellow, minimum size=4pt, inner sep=0pt},
  ]
  \tikzstyle{level 1}=[sibling distance=18mm]
\tikzstyle{level 2}=[sibling distance=12mm]
\tikzstyle{level 3}=[sibling distance=6mm]
  \node [black, label=west:{$0$}] (root) {}
      child { node [black, label=west:{$1$}] {}
        child { node [black, label=west:{$3$}] {}
        }
    }
    child { node [black, label=west:{$2$}] {}
        child { node [black, label=west:{$6$}] {}
        }
    }
    child { node [black, label=west:{$4$}] {}
        child { node [black, label=west:{$7$}] {}
        }
        child { node [black, label=west:{$5$}] {}
        }
    }
    ;
\end{tikzpicture}

In this paper we would like to study vertices by their rank. A vertex has rank $k$ if a minimal downward path to a leaf has $k$ edges. Leaves have rank 0. Surprisingly there is connection between vertices of rank 1 and derangements of size $n$. A \textit{derangement} is a permutation with no fixed points, i.e. a permutation with no cycles of size 1.

In Exercise 3.73 of \cite{bona2022}, it is proved that the number of vertices with rank 1 in increasing trees of size $n$ is equal to the number of derangements of size $n$. However, that proof uses a generating function argument, and a combinatorial proof is certainly desirable for such a surprising result. We provide such a proof.

\section{2. Bijection between vertices of rank 1 and derangements}{Bijection between vertices of rank 1 and derangements}\label{asymptotics}

Let $D_n$ be the set of derrangements of size $n$ and let $T_n$ be the set of increasing trees of size $n$ with a marked vertex of rank 1. We will denote marked vertices as squares colored in red. All permutations will be denoted in disjoint cycle notation.

{\textbf{Theorem.}} For all positive integers $n$, the equality $|D_n|=|T_n|$ holds.

\emph{Proof.} The statement holds trivially for $n=1$. We prove the statement bijectively. For all $n$, we will construct a bijection $f_n:D_n\mapsto A_n$. We will construct $f_n$ recursively. As initial cases we define $f_2$ by 
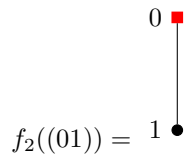
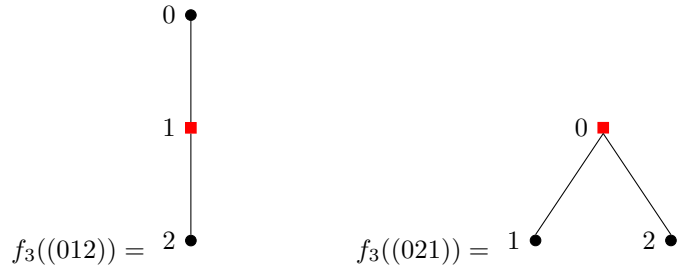
\begin{center}
$f_2((01))=$
\begin{tikzpicture}
  [
    black/.style={circle, draw=black, fill=black, minimum size=4pt, inner sep=0pt},
    red/.style={draw=red, fill=red, minimum size=4pt, inner sep=0pt}
  ]
  \tikzstyle{level 1}=[sibling distance=18mm]
\tikzstyle{level 2}=[sibling distance=12mm]
\tikzstyle{level 3}=[sibling distance=6mm]
  \node [red, label=west:{$0$}] (root) {}
      child { node [black, label=west:{$1$}] {}
    }
    ;
\end{tikzpicture}
\end{center}

and $f_3$ by
\begin{center}
$f_3((012))=$
\begin{tikzpicture}
  [
    black/.style={circle, draw=black, fill=black, minimum size=4pt, inner sep=0pt},
    red/.style={draw=red, fill=red, minimum size=4pt, inner sep=0pt}
  ]
  \tikzstyle{level 1}=[sibling distance=18mm]
\tikzstyle{level 2}=[sibling distance=12mm]
\tikzstyle{level 3}=[sibling distance=6mm]
  \node [black, label=west:{$0$}] (root) {}
      child { node [red, label=west:{$1$}] {}
        child { node [black, label=west:{$2$}] {}
        }
    }
    ;
\end{tikzpicture}\hspace{2cm}  
$f_3((021))=$
\begin{tikzpicture}
  [
   black/.style={circle, draw=black, fill=black, minimum size=4pt, inner sep=0pt},
    red/.style={draw=red, fill=red, minimum size=4pt, inner sep=0pt}
  ]
  \tikzstyle{level 1}=[sibling distance=18mm]
\tikzstyle{level 2}=[sibling distance=12mm]
\tikzstyle{level 3}=[sibling distance=6mm]
  \node [red, label=west:{$0$}] (root) {}
      child { node [black, label=west:{$1$}] {}
    }
    child { node [black, label=west:{$2$}] {}
    }
    ;
\end{tikzpicture}
\end{center}

We split the definiton of $f_n$ into two main cases: if $n-1$ appears in a two-cycle and if $n-1$ appears in a cycle of length greater than two. We will further split these into subcases, each with their own example.

\textbf{Case 1).} The entry $n-1$ is not in a 2-cycle.

If we have $p\in D_n$ such that $n-1$ does not appear in a 2-cycle, then we remove $n-1$ from $p$ written in cycle notation and we are left with $p'\in D_{n-1}$. We now look at the tree $t=f_{n-1}(p')$. Let $k$ be the marked vertex of $t$. We now attach $n-1$ to the vertex $v=p^{-1}(n-1)$ to get a new tree $t'$ with the same marked vertex $k$. We further split this into cases, each with their own example.

\textbf{Case 1)a).} The marked vertex $k=v$, i.e. $n-1$ is now attached to the marked vertex.

Note that in this case, the marked vertex $k=v$ still has rank 1 and it must have at least two leaf children, one of which is $n-1$. Set $f_n(p)=t'$.

\textbf{Example Case 1)a).} Let $p=(053)(142)\in D_6$. The entry $5$ is not in a 2-cycle so we apply the process described in the previous paragraph. Remove 5 to get $p'=(03)(142)$. It has not been proven yet but we will show that 
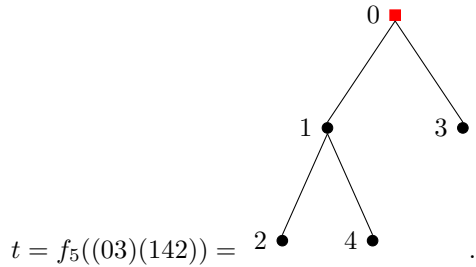
\begin{center}
$t=f_5((03)(142))=$
\begin{tikzpicture}
  [
    black/.style={circle, draw=black, fill=black, minimum size=4pt, inner sep=0pt},
    red/.style={draw=red, fill=red, minimum size=4pt, inner sep=0pt}
  ]
  \tikzstyle{level 1}=[sibling distance=18mm]
\tikzstyle{level 2}=[sibling distance=12mm]
\tikzstyle{level 3}=[sibling distance=6mm]
  \node [red, label=west:{$0$}] (root) {}
      child { node [black, label=west:{$1$}] {}
        child { node [black, label=west:{$2$}] {}
        }  
        child { node [black, label=west:{$4$}] {}
        }  
    }
    child { node [black, label=west:{$3$}] {}
    }
    ;
\end{tikzpicture}.
\end{center}
Then we attach $5$ to $p^{-1}(5)=0$ (in cycle notation it is simply the entry that precedes 5 in its cycle). We get the tree 
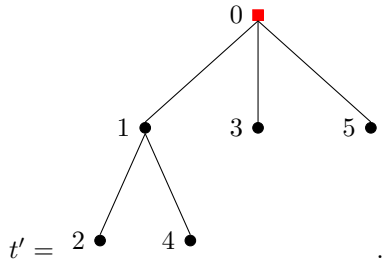
\begin{center}
$t'=$
\begin{tikzpicture}
  [
    black/.style={circle, draw=black, fill=black, minimum size=4pt, inner sep=0pt},
    red/.style={draw=red, fill=red, minimum size=4pt, inner sep=0pt}
  ]
  \tikzstyle{level 1}=[sibling distance=15mm]
\tikzstyle{level 2}=[sibling distance=12mm]
\tikzstyle{level 3}=[sibling distance=6mm]
  \node [red, label=west:{$0$}] (root) {}
      child { node [black, label=west:{$1$}] {}
        child { node [black, label=west:{$2$}] {}
        }  
        child { node [black, label=west:{$4$}] {}
        }  
    }
    child { node [black, label=west:{$3$}] {}
    }
    child { node [black, label=west:{$5$}] {}
        } 
    ;
\end{tikzpicture}.
\end{center}
Finally we set $f_6((053)(142))=t'$.

\textbf{Case 1)b).} The vertex $v$ is neither $k$ nor a child of $k$.

In this case, it is not possible that attaching $n-1$ to $v$ changes the rank of $k$, so $t'$ still has a marked vertex $k$ of rank 1. Set $f_n(p)=t'$.

\textbf{Example Case 1)b).} Let $p=(0523)(14)\in D_6$. The entry $5$ is not in a 2-cycle so we apply the process described in the previous paragraph. Remove 5 to get $p'=(023)(14)$. We can show that 
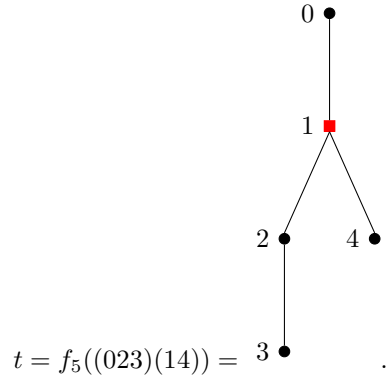
\begin{center}
$t=f_5((023)(14))=$
\begin{tikzpicture}
  [
    black/.style={circle, draw=black, fill=black, minimum size=4pt, inner sep=0pt},
    red/.style={draw=red, fill=red, minimum size=4pt, inner sep=0pt}
  ]
  \tikzstyle{level 1}=[sibling distance=18mm]
\tikzstyle{level 2}=[sibling distance=12mm]
\tikzstyle{level 3}=[sibling distance=6mm]
  \node [black, label=west:{$0$}] (root) {}
      child { node [red, label=west:{$1$}] {}
        child { node [black, label=west:{$2$}] {}
                child { node [black, label=west:{$3$}] {}
                }
        }  
        child { node [black, label=west:{$4$}] {}
        }  
    }
    ;
\end{tikzpicture}.
\end{center}
Then we attach $5$ to $p^{-1}(5)=0$. We get the tree 
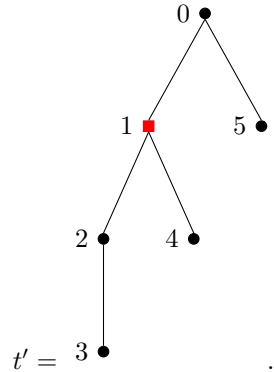
\begin{center}
$t'=$
\begin{tikzpicture}
  [
    black/.style={circle, draw=black, fill=black, minimum size=4pt, inner sep=0pt},
    red/.style={draw=red, fill=red, minimum size=4pt, inner sep=0pt}
  ]
  \tikzstyle{level 1}=[sibling distance=15mm]
\tikzstyle{level 2}=[sibling distance=12mm]
\tikzstyle{level 3}=[sibling distance=6mm]
  \node [black, label=west:{$0$}] (root) {}
      child { node [red, label=west:{$1$}] {}
        child { node [black, label=west:{$2$}] {}
                child { node [black, label=west:{$3$}] {}
                }
        }  
        child { node [black, label=west:{$4$}] {}
        }  
    }
    child { node [black, label=west:{$5$}] {}
        } 
    ;
\end{tikzpicture}.
\end{center}
Finally we set $f_6((0523)(14))=t'$.

\textbf{Case 1)c).} The vertex $v$ is a child of $k$, i.e. $n-1$ is attached to a child of $k$.

We further split this case into two subcases.

\textbf{Case 1)c)i). } The marked vertex $k$ still has rank 1 after attaching $n-1$.

In this case $k$ must have a leaf child which is not $v$. The marked vertex $k$ is still has rank 1, so set $f_n(p)=t'$.

\textbf{Example Case 1)c)i).} Let $p=(031425)\in D_6$. The entry $5$ is not in a 2-cycle. Remove 5 to get $p'=(03142)$. It has not been proven yet but we will show that 
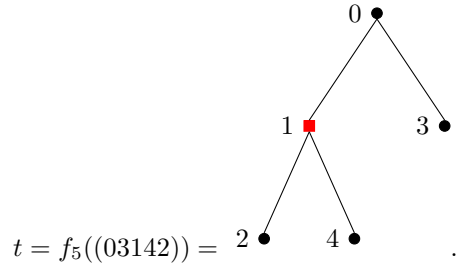
\begin{center}
$t=f_5((03142))=$
\begin{tikzpicture}
  [
    black/.style={circle, draw=black, fill=black, minimum size=4pt, inner sep=0pt},
    red/.style={draw=red, fill=red, minimum size=4pt, inner sep=0pt}
  ]
  \tikzstyle{level 1}=[sibling distance=18mm]
\tikzstyle{level 2}=[sibling distance=12mm]
\tikzstyle{level 3}=[sibling distance=6mm]
  \node [black, label=west:{$0$}] (root) {}
      child { node [red, label=west:{$1$}] {}
        child { node [black, label=west:{$2$}] {}
        }  
        child { node [black, label=west:{$4$}] {}
            }  
        }
    child { node [black, label=west:{$3$}] {}
                }
    ;
\end{tikzpicture}.
\end{center}
Then we attach $5$ to $p^{-1}(5)=2$. We get the tree 
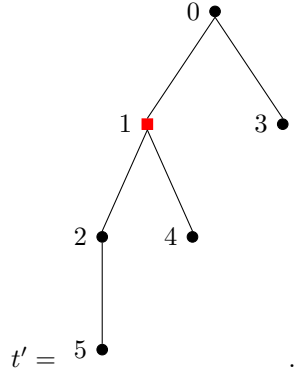
\begin{center}
$t'=$
\begin{tikzpicture}
  [
    black/.style={circle, draw=black, fill=black, minimum size=4pt, inner sep=0pt},
    red/.style={draw=red, fill=red, minimum size=4pt, inner sep=0pt}
  ]
  \tikzstyle{level 1}=[sibling distance=18mm]
\tikzstyle{level 2}=[sibling distance=12mm]
\tikzstyle{level 3}=[sibling distance=6mm]
  \node [black, label=west:{$0$}] (root) {}
      child { node [red, label=west:{$1$}] {}
        child { node [black, label=west:{$2$}] {}
            child { node [black, label=west:{$5$}] {}
            }
        }  
        child { node [black, label=west:{$4$}] {}
            }  
        }
    child { node [black, label=west:{$3$}] {}
                }
    ;
\end{tikzpicture}.
\end{center}
Finally we set $f_6((031425))=t'$.

\textbf{Case 1)c)ii).} The marked vertex $k$ does not have rank 1 after attaching $n-1$.

This means $v$ was the only leaf child of $k$. In $t'$, $k$ has rank 2 and $v$ has rank $1$. Transform $t'$ into $t''$ by marking $v$. Set $f_n(p)=t''$. Note that in this case $n-1$ is attached to the marked vertex, $n-1$ is the only child of the marked vertex, and the parent of the marked vertex has rank 2.

\textbf{Example Case 1)c)ii).} Let $p=(013524)\in D_6$. The entry $5$ is not in a 2-cycle. Remove 5 to get $p'=(01324)$. We can show that 
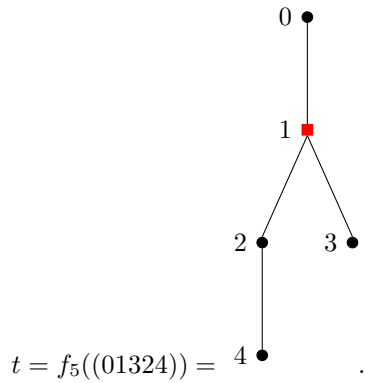
\begin{center}
$t=f_5((01324))=$
\begin{tikzpicture}
  [
    black/.style={circle, draw=black, fill=black, minimum size=4pt, inner sep=0pt},
    red/.style={draw=red, fill=red, minimum size=4pt, inner sep=0pt}
  ]
  \tikzstyle{level 1}=[sibling distance=18mm]
\tikzstyle{level 2}=[sibling distance=12mm]
\tikzstyle{level 3}=[sibling distance=6mm]
  \node [black, label=west:{$0$}] (root) {}
      child { node [red, label=west:{$1$}] {}
        child { node [black, label=west:{$2$}] {}
            child { node [black, label=west:{$4$}] {}
                }
        }  
        child { node [black, label=west:{$3$}] {}
            }  
        }
    ;
\end{tikzpicture}.
\end{center}
Then we attach $5$ to $p^{-1}(5)=3$. We get the tree 
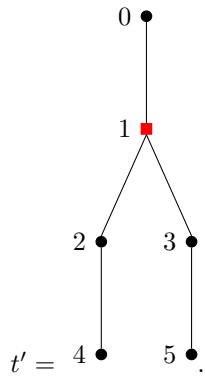
\begin{center}
$t'=$
\begin{tikzpicture}
  [
    black/.style={circle, draw=black, fill=black, minimum size=4pt, inner sep=0pt},
    red/.style={draw=red, fill=red, minimum size=4pt, inner sep=0pt}
  ]
  \tikzstyle{level 1}=[sibling distance=18mm]
\tikzstyle{level 2}=[sibling distance=12mm]
\tikzstyle{level 3}=[sibling distance=6mm]
  \node [black, label=west:{$0$}] (root) {}
      child { node [red, label=west:{$1$}] {}
        child { node [black, label=west:{$2$}] {}
            child { node [black, label=west:{$4$}] {}
                }
        }  
        child { node [black, label=west:{$3$}] {}
                child { node [black, label=west:{$5$}] {}
            }
            }  
        }
    ;
\end{tikzpicture}.
\end{center}
The marked vertex has rank 2, so we instead mark the parent of $n-1=5$ to get $t''$.
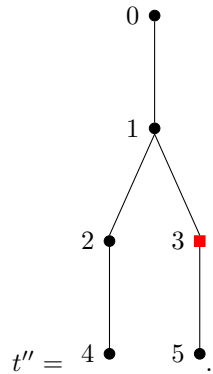
\begin{center}
$t''=$
\begin{tikzpicture}
  [
    black/.style={circle, draw=black, fill=black, minimum size=4pt, inner sep=0pt},
    red/.style={draw=red, fill=red, minimum size=4pt, inner sep=0pt}
  ]
  \tikzstyle{level 1}=[sibling distance=18mm]
\tikzstyle{level 2}=[sibling distance=12mm]
\tikzstyle{level 3}=[sibling distance=6mm]
  \node [black, label=west:{$0$}] (root) {}
      child { node [black, label=west:{$1$}] {}
        child { node [black, label=west:{$2$}] {}
            child { node [black, label=west:{$4$}] {}
                }
        }  
        child { node [red, label=west:{$3$}] {}
                child { node [black, label=west:{$5$}] {}
            }
            }  
        }
    ;
\end{tikzpicture}.
\end{center}
Finally we set $f_6((013524))=t'$.

All of Cases 1)a), 1)b), and 1)b)i) are 'nice' in the sense that removing $n-1$ leaves us with an increasing tree with a marked vertex of rank 1. Case 1)c)ii) is 'almost nice' in the sense that removing $n-1$ leaves us with a marked leaf and the removal of $n-1$ is what caused the parent of the marked leaf to have rank 1.

Now we move on to the next case.

\textbf{Case 2). }The entry $n-1$ in $p\in D_n$ is in a 2-cycle with some entry $j$. 

We remove the 2-cycle $(n-1\text{ }j)$ and we get a derangement $p'$ on the set $\{0,...,n-2\}-\{v\}$. Up to relabeling then $p'$ is a derangement of size $n-2$, so $t=f_{n-2}(p')$ is defined by our recursion. Let $k$ be the marked vertex of $t$. We split Case 2) into two subcases. 

\textbf{Case 2)a). }The marked vertex $k>j$. 

We define the following algorithm which will produce a new tree $t'$ on $n$ vertices.
\begin{enumerate}
    \item Consider the upward path from $k$ to the root in $t$. Attach $j$ to this path in its only possible position given by the increasing characteristic of $t$. If $j=0$, then $j$ will simply be the new root whose only child is the root of $t$. At this point $j$ should only have one child.
    \item Our objective is that a depth first walk starting at $j$ will be such that $k$ will be the first vertex of rank 1 we encounter. Starting at $j$, go down the tree through the downward path to $k$ until you get to a vertex $w$ with rank 1 or with rank 2 or higher with two or more children. If $w$ has rank 1 and $k\neq w$ then it must be that $w$ has other children, one of which is the root of the subtree containing $k$. Remove the subtree containing $k$ and attach it so that the root is a child of $j$. If $w$ is such that it has rank at least 2, look at its children $c_1,...,c_m$. Let $k$ be in the subtree rooted at $c_i$. If $c_i=max\{c_1,...,c_m\}$ leave the tree as it is. Otherwise attach the subtree rooted at $c_i$ to $j$ so that $j$ is the parent of $c_i$. Go down the path.
    \item Repeat Step 2 until you reach $k$.
    \item Attach $n-1$ to $j$ and mark $j$
\end{enumerate}
This results in a tree $t'$ on $n$ vertices with a marked vertex of rank 1. Set $f_n(p)=t'$. Note that in this case $n-1$ is attached to the marked vertex. Moreover we claim that $n-1$ is not an only child and it is the only leaf sibling of its parent. Since $j$ was attached along the path from the root to $k$ it had at least one child when it was attached. The algorithm only increases the number of children of $j$, so $n-1$ is not an only child. Also note that all of the subtrees that were attached to $j$ were so that the root had rank at least 1.

\textbf{Example Case 2)a).} Let $p=(09)(135268)(47)\in D_{10}$. Then $9=n-1$ is in a two-cycle, so remove $(09)$ and get $p'=(135268)(47)$. We can show that 

$f_8(p')=t=$
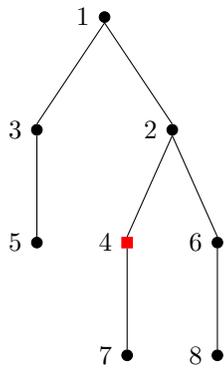
\begin{tikzpicture}
  [
    black/.style={circle, draw=black, fill=black, minimum size=4pt, inner sep=0pt},
    red/.style={draw=red, fill=red, minimum size=4pt, inner sep=0pt}
  ]
  \tikzstyle{level 1}=[sibling distance=18mm]
\tikzstyle{level 2}=[sibling distance=12mm]
\tikzstyle{level 3}=[sibling distance=8mm]
  \node [black, label=west:{$1$}] (root) {}
      child { node [black, label=west:{$3$}] {}
      child { node [black, label=west:{$5$}] {}
        } 
    }
    child { node [black, label=west:{$2$}] {}
        child { node [red, label=west:{$4$}] {}
            child { node [black, label=west:{$7$}] {}
            }
        }
        child { node [black, label=west:{$6$}] {}
            child { node [black, label=west:{$8$}] {}
            }
        }
    }
    ;
\end{tikzpicture}

Then the marked vertex $k=4$ and $j=0$, so $j<k$. The downward path between the marked vertex and root is marked by blue vertices below.

$t=$
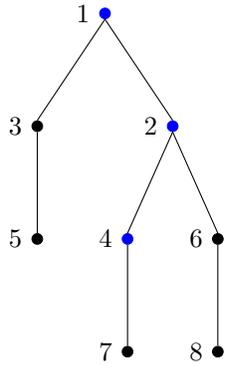
\begin{tikzpicture}
  [
    black/.style={circle, draw=black, fill=black, minimum size=4pt, inner sep=0pt},
    red/.style={draw=red, fill=red, minimum size=4pt, inner sep=0pt},
    blue/.style={circle, draw=blue, fill=blue, minimum size=4pt, inner sep=0pt},
  ]
  \tikzstyle{level 1}=[sibling distance=18mm]
\tikzstyle{level 2}=[sibling distance=12mm]
\tikzstyle{level 3}=[sibling distance=8mm]
  \node [blue, label=west:{$1$}] (root) {}
      child { node [black, label=west:{$3$}] {}
      child { node [black, label=west:{$5$}] {}
        } 
    }
    child { node [blue, label=west:{$2$}] {}
        child { node [blue, label=west:{$4$}] {}
            child { node [black, label=west:{$7$}] {}
            }
        }
        child { node [black, label=west:{$6$}] {}
            child { node [black, label=west:{$8$}] {}
            }
        }
    }
    ;
\end{tikzpicture}

The vertex 0 can only be inserted in this path as the parent of 1.

$t=$
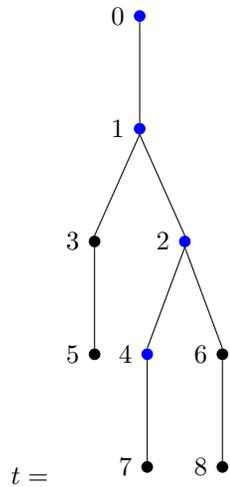
\begin{tikzpicture}
  [
    black/.style={circle, draw=black, fill=black, minimum size=4pt, inner sep=0pt},
    red/.style={draw=red, fill=red, minimum size=4pt, inner sep=0pt},
    blue/.style={circle, draw=blue, fill=blue, minimum size=4pt, inner sep=0pt},
  ]
  \tikzstyle{level 1}=[sibling distance=18mm]
\tikzstyle{level 2}=[sibling distance=12mm]
\tikzstyle{level 3}=[sibling distance=10mm]
  \node [blue, label=west:{$0$}] (root) {}
  child { node [blue, label=west:{$1$}] {}
      child { node [black, label=west:{$3$}] {}
      child { node [black, label=west:{$5$}] {}
        } 
    }
    child { node [blue, label=west:{$2$}] {}
        child { node [blue, label=west:{$4$}] {}
            child { node [black, label=west:{$7$}] {}
            }
        }
        child { node [black, label=west:{$6$}] {}
            child { node [black, label=west:{$8$}] {}
            }
        }
    }
    }
    ;
\end{tikzpicture}

Now we start at 0 and we want to get down to 4. We drop down to 1, where we have a choice between 3 and 2. The marked vertex 4 is in the subtree of 2, but two is less than 3, so we attach 2 and its subtree to 0.

$t=$
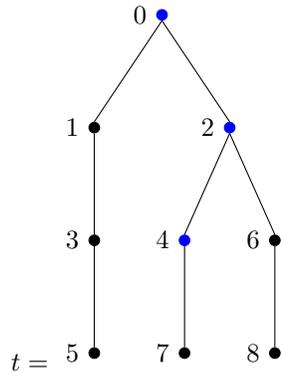
\begin{tikzpicture}
  [
    black/.style={circle, draw=black, fill=black, minimum size=4pt, inner sep=0pt},
    red/.style={draw=red, fill=red, minimum size=4pt, inner sep=0pt},
    blue/.style={circle, draw=blue, fill=blue, minimum size=4pt, inner sep=0pt},
  ]
  \tikzstyle{level 1}=[sibling distance=18mm]
\tikzstyle{level 2}=[sibling distance=12mm]
\tikzstyle{level 3}=[sibling distance=10mm]
  \node [blue, label=west:{$0$}] (root) {}
  child { node [black, label=west:{$1$}] {}
      child { node [black, label=west:{$3$}] {}
      child { node [black, label=west:{$5$}] {}
        } 
    }
    }
   child { node [blue, label=west:{$2$}] {}
        child { node [blue, label=west:{$4$}] {}
            child { node [black, label=west:{$7$}] {}
            }
        }
        child { node [black, label=west:{$6$}] {}
            child { node [black, label=west:{$8$}] {}
            }
        }
    }
    ;
\end{tikzpicture}

Now we go down to 2. Again we have a choice between 4 and 6. We want to go down to 4, but 4 is less than 6, so we attach 4 and its subtree to 0.

$t=$
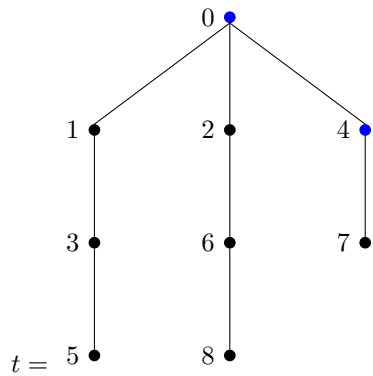
\begin{tikzpicture}
  [
  black/.style={circle, draw=black, fill=black, minimum size=4pt, inner sep=0pt},
    red/.style={draw=red, fill=red, minimum size=4pt, inner sep=0pt},
    blue/.style={circle, draw=blue, fill=blue, minimum size=4pt, inner sep=0pt},
  ]
  \tikzstyle{level 1}=[sibling distance=18mm]
\tikzstyle{level 2}=[sibling distance=12mm]
\tikzstyle{level 3}=[sibling distance=10mm]
  \node [blue, label=west:{$0$}] (root) {}
  child { node [black, label=west:{$1$}] {}
      child { node [black, label=west:{$3$}] {}
      child { node [black, label=west:{$5$}] {}
        } 
    }
    }
   child { node [black, label=west:{$2$}] {}
        child { node [black, label=west:{$6$}] {}
            child { node [black, label=west:{$8$}] {}
            }
        }
    }
    child { node [blue, label=west:{$4$}] {}
            child { node [black, label=west:{$7$}] {}
            }
        }
    ;
\end{tikzpicture}

Next we go down to $4=k$, so we stop. Now we attach 9 to 0 and mark 0.

$t'=$
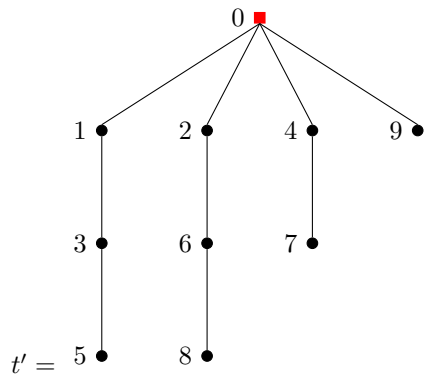
\begin{tikzpicture}
  [
    black/.style={circle, draw=black, fill=black, minimum size=4pt, inner sep=0pt},
    red/.style={draw=red, fill=red, minimum size=4pt, inner sep=0pt},
    blue/.style={circle, draw=blue, fill=blue, minimum size=4pt, inner sep=0pt},
  ]
  \tikzstyle{level 1}=[sibling distance=14mm]
\tikzstyle{level 2}=[sibling distance=12mm]
\tikzstyle{level 3}=[sibling distance=10mm]
  \node [red, label=west:{$0$}] (root) {}
  child { node [black, label=west:{$1$}] {}
      child { node [black, label=west:{$3$}] {}
      child { node [black, label=west:{$5$}] {}
        } 
    }
    }
   child { node [black, label=west:{$2$}] {}
        child { node [black, label=west:{$6$}] {}
            child { node [black, label=west:{$8$}] {}
            }
        }
    }
    child { node [black, label=west:{$4$}] {}
            child { node [black, label=west:{$7$}] {}
            }
        }
    child { node [black, label=west:{$9$}] {}
    }
    ;
\end{tikzpicture}

So $f_{10}((09)(135268)(47))=t'$

\textbf{Case 2)b). }The marked vertex $k<j$.

This case is much simpler. Attach $j$ to $k$ and then attach $n-1$ to $j$ and mark $j$ to get $t'$. Set $f_n(p)=t'$. Note that in this case $n-1$ is attached to the marked vertex. Moreover $n-1$ is an only child and the parent of the marked vertex has rank 1. Indeed, the parent of the marked vertex is $k$ which was marked in $t$, so it had rank 1.

\textbf{Example Case 2)b). } Let $p=(025168)(37)(49)\in D_{10}$. Then $9=n-1$ is in a two-cycle, so remove $(49)$ and get $p'=(025168)(37)$. We can show that 

$f_8(p')=t=$
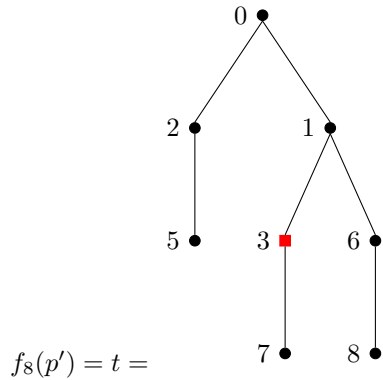
\begin{tikzpicture}
  [
    black/.style={circle, draw=black, fill=black, minimum size=4pt, inner sep=0pt},
    red/.style={draw=red, fill=red, minimum size=4pt, inner sep=0pt},
    blue/.style={circle, draw=blue, fill=blue, minimum size=4pt, inner sep=0pt},
  ]
  \tikzstyle{level 1}=[sibling distance=18mm]
\tikzstyle{level 2}=[sibling distance=12mm]
\tikzstyle{level 3}=[sibling distance=8mm]
  \node [black, label=west:{$0$}] (root) {}
      child { node [black, label=west:{$2$}] {}
      child { node [black, label=west:{$5$}] {}
        } 
    }
    child { node [black, label=west:{$1$}] {}
        child { node [red, label=west:{$3$}] {}
            child { node [black, label=west:{$7$}] {}
            }
        }
        child { node [black, label=west:{$6$}] {}
            child { node [black, label=west:{$8$}] {}
            }
        }
    }
    ;
\end{tikzpicture}

In $t$, $k=3<4=j$. So we attach $4$ to 3 and attach 9 to 4 to get $t'$

$t'=$
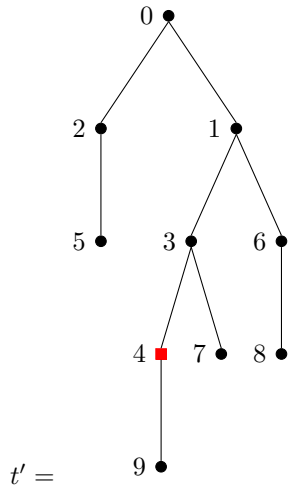
\begin{tikzpicture}
  [
    black/.style={circle, draw=black, fill=black, minimum size=4pt, inner sep=0pt},
    red/.style={draw=red, fill=red, minimum size=4pt, inner sep=0pt},
    blue/.style={circle, draw=blue, fill=blue, minimum size=4pt, inner sep=0pt},
  ]
  \tikzstyle{level 1}=[sibling distance=18mm]
\tikzstyle{level 2}=[sibling distance=12mm]
\tikzstyle{level 3}=[sibling distance=8mm]
  \node [black, label=west:{$0$}] (root) {}
      child { node [black, label=west:{$2$}] {}
      child { node [black, label=west:{$5$}] {}
        } 
    }
    child { node [black, label=west:{$1$}] {}
        child { node [black, label=west:{$3$}] {}
        child { node [red, label=west:{$4$}] {}
                child { node [black, label=west:{$9$}] {}
            }
            }
            child { node [black, label=west:{$7$}] {}
            }
        }
        child { node [black, label=west:{$6$}] {}
            child { node [black, label=west:{$8$}] {}
            }
        }
    }
    ;
\end{tikzpicture}

The following diagram shows all the possibilities of trees in the image of $f_n$. The parent of $n-1$ is denoted $v$ and the marked vertex is $k$.
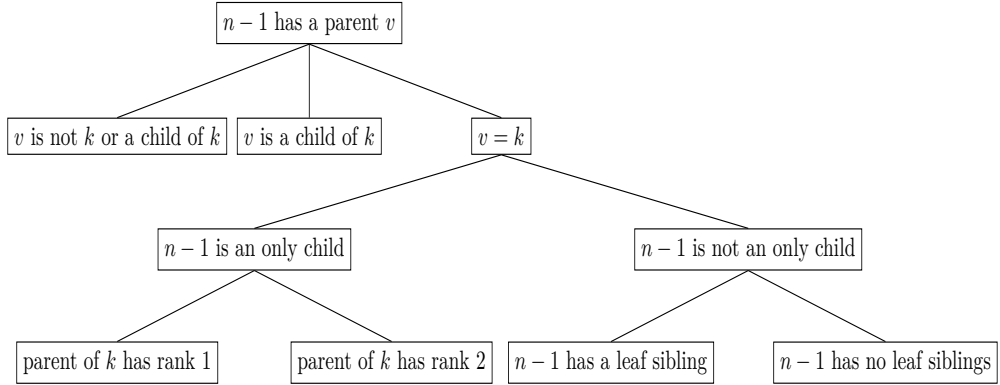
\begin{center}
\resizebox{5.25in}{2in}{
\begin{tikzpicture}[sibling distance=10em,
  every node/.style = {shape=rectangle,
    draw, align=center}]
    \tikzstyle{level 1}=[sibling distance=35mm]
    \tikzstyle{level 2}=[sibling distance=90mm]
    \tikzstyle{level 3}=[sibling distance=50mm]
  \node {$n-1$ has a parent $v$}
    child { node {$v$ is not $k$ or a child of $k$} }
    child { node {$v$ is a child of $k$} }
    child { node {$v=k$} 
        child{ node{$n-1$ is an only child}
            child{ node{parent of $k$ has rank 1}
            } 
            child{ node{parent of $k$ has rank 2}
            }
        }
        child{ node{$n-1$ is not an only child}
            child{ node{$n-1$ has a leaf sibling}
            } 
            child{ node{$n-1$ has no leaf siblings}
            }
        }
    }
    ;
\end{tikzpicture}
}
\end{center}
Note that every increasing tree with a marked vertex of rank 1 is covered by exactly one of the cases on the outer nodes in the diagram above. Moreover, each outer node covers one case from the construction of our function $f_n$, as seen in the diagram below.
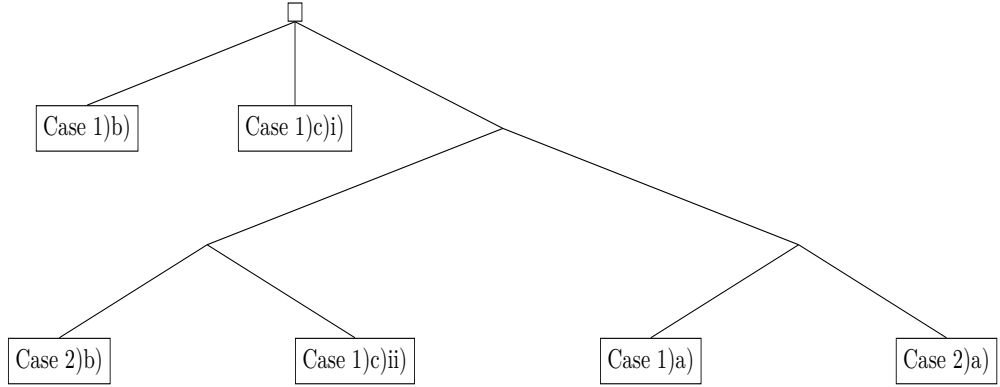
\begin{center}
\resizebox{5.25in}{2in}{
\begin{tikzpicture}
[sibling distance=10em,
  every node/.style = {shape=rectangle,
    draw, align=center}]]
    \tikzstyle{level 2}=[sibling distance=100mm]
    \tikzstyle{level 3}=[sibling distance=50mm]
  \node {}
    child { node {Case 1)b)} }
    child { node {Case 1)c)i)} }
    child { 
        child{ 
            child{ node{Case 2)b)}
            } 
            child{ node{Case 1)c)ii)}
            }
        }
        child{ 
            child{ node{Case 1)a)}
            } 
            child{ node{Case 2)a)}
            }
        }
    }
    ;
\end{tikzpicture}
}
\end{center}
Let $A_n$ be the set of derangements of size $n$ in which $n-1$ is not in a 2-cycle and let $B_n$ be the set of derangements of size $n$ in which $n-1$ is in a 2-cycle. 

We first show that $f_n$ restricted to $A_n$ is a bijection. Consider the trees obtained from Case 1)a) through c), which correspond to trees in the image of $A_n$. For any tree in Case 1)a), Case 1)b), and Case 1)c)i) we can remove $n-1$ and get a tree $t$ with a marked vertex of rank 1. From this tree we get the derangement $f_{n-1}^{-1}(t)$ and and we insert $n-1$ after its parent in $t$ in the cycles of $f_{n-1}^{-1}(t)$. This process shows how to recover the preimage of any such marked tree in a unique way, so $f_n$ is bijective when restricted to these cases.

For the Case 1)c)ii), in $t$, $n-1$ is attached to the marked vertex $k$, and its parent $p$ has rank 2. So when we remove $n-1$ to get $t'$, the marked vertex $k$ is now a leaf and $p$ has rank 1. Moreover, $k$ is the only leaf child of $p$. We get $f_{n-1}^{-1}(t')$ and insert $n-1$ after its parent in $t$ in the cycles of $f_{n-1}^{-1}(t)$. So $f_n$ is bijective when restricted to $A_n$.

To show $f_n$ is bijective when restricted to $B_n$, consider first trees obtained from Case 2)b). In this case, $n-1$ is attached to the marked vertex $k$, $n-1$ is an only child, and the parent $p$ of $k$ has rank 1. Then if we remove $n-1$ and $k$ and mark $p$ we get a marked tree $t$ on $n-2$ vertices. We can get the derangement $f_{n-2}^{-1}(t)$ and add the two-cycle with $n-1$ and $k$ to get a derangement of size $n$. This process shows how to recover the preimage of any such marked tree in a unique way, so $f_n$ is bijective when restricted to Case 1)b).

Finally we look at trees obtained from Case 2)a). In these trees, $n-1$ is attached to the marked vertex $k$ and it has siblings, none of which are leaves. Let $c_1,...,c_m$ be the nonleaf children of $k$ ordered from least to greatest. To recover the preimage of such a tree, we perform a depth first walk starting at $k$ until we reach a vertex of rank 1. 

Once we stop, we want to join all the subtrees rooted at $c_1,...,c_m$ and bring them as one single tree with root $c_1$. To do this, start at $c_{m-1}$ and go down a first depth walk until you reach a vertex $q$ which has rank 1 or has rank at least 2 such that one of its children is greater than $c_m$. Then attach $c_m$, and its subtree, as a child of $q$. Repeat the process until you attach $c_2$ somewhere in the subtree rooted at $c_1$. Then remove $k$ and $n-1$ from the tree and have $c_1$ take the place of $k$ in the new tree. Then we obtain a marked tree on $n-2$ vertices. We recover the preimage in $f_{n-2}$ and attach the 2-cycle with $n-1$ and $k$ in it. The process described above reverses the algorithm from Case 2)b) because we attach the vertices $c_1,...,c_m$ in a way that taking the first depth walk in the algorithm would have detached them. So we recover the preimage in a unique way, and so $f_n$ is bijective when restricted to $B_n$. Since $D_n=A_n\sqcup B_n$, $f_n$ is a bijection. \qed

\section{3. Conclusion}\label{Conclusion}

The result presented above was originally obtained by way of solving a differential equation to obtain the exponential generating function for vertices of rank 1. The same method can be used to obtain a differential equation with the exponential generating function for vertices of rank 2, as discussed verbally by the author with B{\'o}na. Through this method it has been shown that the generating function can not be expressed in terms of elementary functions. However, it is very likely that the vertices of rank 2 are in correspondence to some other statistic or class of permutations. What we have shown in this paper is essentially that the vertices of rank 1 satisfy the same recurrence as derangements, $A_n=nA_{n-1}+nA_{n-2}$. A further direction of research would be to formulate a similar recurrence for vertices of rank 2 and, possibly, general rank $k$.

It would also be of interest to count trees (or permutations) corresponding to the six individual cases in the proof. For instance, we could count the number of vertices of rank 1 that have $n-1$ as a child by finding the total number of trees corresponding to the cases 2)a), 2)b), 1)a), and 1)c)ii). Out of these cases only 1)c)ii) seems to be nontrivial.

\printbibliography

\end{document}